\documentclass{elsart}
\usepackage{amssymb,amsmath,epsfig}
\newcommand{\tsm}[1]{\left(\begin{smallmatrix}#1\end{smallmatrix}\right)}
\let\MySum=\sum
\def\sum[#1,#2]{\MySum_{#1}^{#2}}
\newcommand{\ma}[1]{\displaystyle\left(\begin{matrix}#1\end{matrix}\right)}
\newcommand{\tma}[1]{\textstyle\left(\begin{matrix}#1\end{matrix}\right)}
\newcommand{\sm}[1]{\displaystyle\begin{matrix}#1\end{matrix}}
\journal{Applied Numerical Mathematics}
\volume{42}
\pubyear{2002}
\newcounter{mylastpage}
\issue{1--3}
\usepackage{hyperref}
\makeatletter
\def\ps@copyright{%
 \let\@oddhead\@empty
 \let\@evenhead\@empty
 \def\@oddfoot{\small\slshape\hskip-7em
   Published in \@journal\ \@volume\ (\the\@pubyear)\ no.\ \@issue, pp.\ \ESpagenumber{firstpage}--\ESpagenumber{mylastpage},
\href{http://dx.doi.org/10.1016/S0168-9274(01)00157-X}{doi: 10.1016/S0168-9274(01)00157-X}}%
 \let\@evenfoot\@oddfoot}
\makeatother
\begin{document}

\begin{frontmatter}

\title{On quasi--linear PDAEs with convection:\\[2ex]
applications, indices, numerical solution}
\author{W. Lucht},
\ead{lucht@mathematik.uni-halle.de}
\author{K. Debrabant}
\address{ Martin-Luther-Universit\"{a}t Halle--Wittenberg,\\
  Fachbereich Mathematik und Informatik,
  Institut f\"{u}r Numerische Mathematik,\\
   Postfach, D--06099 Halle, Germany }

\begin{abstract}
For a class of partial differential algebraic equations (PDAEs)
of quasi--linear type  which include nonlinear terms of
convection type a possibility to determine a time and spatial
index is considered. As a typical example we investigate an application
from plasma physics. Especially
we discuss the numerical solution of
initial boundary value problems by means of a corresponding
finite difference splitting procedure which is a modification of
a well known  fractional step method coupled with a matrix
factorization. The convergence of the numerical solution towards
the exact solution of the corresponding initial boundary value
problem is investigated. Some  results of a numerical solution of
the plasma PDAE are given.
\end{abstract}

\end{frontmatter}

\section{Introduction} \label{introd}
In this paper quasi--linear PDAEs for $u=u(t,x),~x\in \Omega:=(0,1),$ of the form
\begin{align}\label{pdae1}
Au_t+Bu_{xx}+C[u]u_x+Du=&f(t,x),\quad t\in (0,t_e),
~~~~x\in\Omega,
\end{align}
for some $t_e>0$ are considered.
$~u$ and $f$ are mappings
$~u,f:~[0,t_e]\times\bar\Omega\rightarrow{\Rset}^n,~~n\geq 1,$ where  $f$ (supposed to be
sufficiently smooth) is given.
$A,~B,~C[u]$ and $D$ are real $(n,n)-$matrices where $A,~B$ and $D$
are assumed to be constant.  All matrices may be
singular, but $A,~B\neq 0$.
$C[u]$ may depend on $u$.
We suppose that it is linear in $u$. Typically, when $C[u]$  is linear  in $u$, vector
$C[u]u_{x}$ describes physical convection.\\[1ex]
For system (\ref{pdae1}) we study classical initial boundary value problems (IBVPs).
Initial values (IVs) may be decomposed
\begin{align}\label{IC}
u(0,x)&=\Phi_a(x)+\Phi_c(x),\quad x\in \bar \Omega\\
\nonumber
\Phi_{a,k}(x)=&\begin{cases} u_k(0,x)\quad u_k(0,x) \quad
\hbox{can be prescribed arbitrarily}\\0\hspace{0.8cm}\qquad
\hbox{otherwise,}\end{cases}\\\nonumber
\Phi_{c,k}(x)=&\begin{cases} u_k(0,x)\quad u_k(0,x) \quad \hbox{cannot be
prescribed arbitrarily}\\0\hspace{0.8cm}\qquad \hbox{otherwise,}\end{cases}
\end{align}
with ($k=1,...,n$).
The boundary values (BVs) are of  similar form,
\begin{align}\label{BC} u(t,x)=&\Psi_a(t,x)+\Psi_c(t,x),\quad t\in
[0,t_e],\quad x\in\partial\Omega=\{0,1\}.
\end{align}
This means the data which can be prescribed arbitrarily are in
$\Phi_a,~\Psi_a$. The consistent data (see e.g. \cite{Camp95a} or \cite{Lucht99})
are collected in $\Phi_c,~\Psi_c$.
Furthermore, we assume  that the compatibility relations
\begin{align}\label{com}
\Phi_a(x)+ \Phi_c(x)=\Psi_a(0,x)+\Psi_c(0,x),\quad x\in\partial\Omega,
\end{align}
are satisfied. Because almost every PDAE has its own IV and BV
distribution we avoid here to give a more detailed general description of
these data. In section \ref{numex} we solve this problem for the plasma
PDAE considered in the next section.\\
This paper is organized as
follows. In section \ref{plasma} a nonlinear PDAE
from physics is presented. In section \ref{indices}  known
general index concepts applicable also to nonlinear PDAEs are considered.
The numerical solution of IBVPs by a finite difference splitting method is
studied in section \ref{splitting}. Especially, convergence results are
given. A numerical example from plasma physics is presented in section
\ref{numex}.
\section{An application}\label{plasma}
In this section we cite a mathematical model from plasma physics \cite[Chapter 5]{Dodd} whose
underlying system of equations is of type (\ref{pdae1}).
It describes the space--time--movement of a system consisting of ions (with
mass $m_i$, density $n_i$ and positive electrical charge $q$) and electrons
(with density $n_e$ and electrical
charge $-q$) in the space domain $\Omega$. The charge $q(n_e-n_i)$ produces
an electrical potential $\phi$ such  that the ions move under the
corresponding electrical force $-q\phi_x$. The system of electrons
is considered to be a gas with (constant) temperature $T_e$ and pressure
$p=k_BT_en_e$ where $k_B$ is Boltzmann's constant. The relation between
$n_e$ and $\phi$ is given by the equilibrium of electrical and mechanical forces,
\[qn_e\phi_x-k_BT_en_{e,x}=0,~\]
and the equations of conservation of mass and momentum for the ions are
\[n_{i,t}+(n_iv_i)_x-D_in_{i,xx}=0\quad\mathrm{and}\quad
m_i\left(\frac{\partial}{\partial t}+v_i\frac{\partial}{\partial
x}\right)v_i+q\phi_x=0,\] respectively.
$D_i$ is a (constant) diffusion coefficient, and $v_i$ is the
velocity of the ions. The equation for $\phi$ is Poisson's equation
\[\phi_{xx}-4\pi(n_e-n_i)=0.\]
It is convenient to transform  this system of equations by a linear transformation
into a new system of the form (\ref{pdae1}) with $n=4,~~f=0$ and new dependent variables
$u_i$ and with matrices
\begin{align} \label{matrices}
A=\tsm {1 & 0 & 0 & 0 \\0 & 1 & 0 & 0\\0 & 0 & 0 & 0 \\0&0&0&0},~
B=&\tsm{-b_0&0&0&0\\0&0&0&0\\0&0&0&0\\0&0&0&1},~
C[u]=\tsm{u_2&u_1&0&0\\0&u_2&0&d_1\\0&0&-1&u_3\\0&0&0&0},~
D=\tsm{0&0&0&0\\0&0&0&0\\0&0&0&0\\1&0&-1&0}.
\end{align}
The new variables $u_i$ are (up to a constant) given
by $u_1\sim n_i,~~u_2\sim v_i,~~u_3\sim n_e,~~u_4\sim \phi.$
$b_0\geq 0$ and $d_1>0$ are  constants.\\
Note that the matrices $A,~B$ are diagonal and singular, and $C[u]$ is
linear in the components of $u$ (for short we say that $C[u]$ is linear in
$u$).\\
Further examples of a system of type (\ref{pdae1}) are a poroelastic model of a living bone
\cite{Cowin99,Det93,Lucht001} and the incompressible Navier--Stokes equations
if  eq. (\ref{pdae1}) is generalized in obvious manner to
two and three space dimensions.
\section{Indices}\label{indices}
First we introduce two definitions of classical linear spaces. Let $l,~m$ be nonnegative
integers, and let $C_{sc}^{l,m}$ be the space
of scalar real functions $w(t,x),~~ t\in[0,t_e],~~x\in\overline\Omega$,
whose derivatives up to the $(l+m)^{th}$ order (time derivatives up to the $l^{th}$
order and space derivatives up to the $m^{th}$ order) are continuous.
With $C_{sc}^{l,m}$ the set
\begin{align*}
C^{l,m}_n:=\left\{w=(w_1,...,w_n)^T;~~~w_i\in C_{sc}^{l,m},
~~i=1,...,n\right\}
\end{align*} is defined.
By $C_{n,0}^{l,m}$ we denote a set of vector valued functions with
vanishing BVs,
\begin{align*}
 C^{l,m}_{n,0}:=
\Big\{w\in C_n^{l,m}; \quad &  w(t,0)=w(t,1)=0\Big\}.
\end{align*}
While the index of linear PDAEs has been considered by
several authors, e.g.  \cite{Camp95a,Camp96,Camp96a,Guen98,Lucht99,Mar97},
the index for nonlinear
systems is investigated little, see however \cite{Lucht98} where a
discretization based index definition has been given, and \cite{Martin98}.
We mention that we do not transform system (\ref{pdae1}) to a
first order system because we prefer in numerical calculations the original
second order form.\\ In this paper, a well known index concept for PDAEs
is used to construct certain operators whose invertibility (if so) yields
indices. The basic notion of a  time index $\nu_{t}$ and a spatial index
$\nu_{x}$ can be found, e.g., in \cite{Camp96,Lucht99,Martin98}.\\
Throughout this paper the time index $\nu_t$ of (\ref{pdae1}) is of special interest.
We assume that a solution u of the
IBVP (\ref{pdae1}) -- (\ref{com}) exists and that $u$ is
sufficiently differentiable. Furthermore, we suppose
that $u(t,x)=0~$ for
$~x\in\partial \Omega.$ When $\Psi_a$ and $\Psi_c$ (see eq. (\ref{BC})) are
known, zero BVs can be obtained by a suitable transformation of $u$.
\subsection{Time index}\label{timeindex}
\begin{defn} \label{deftime}
{\it If the matrix $A$ is regular, the  time index $\nu_t$ of the
PDAE (\ref{pdae1}) is defined to be zero. If $A$ is singular, then $\nu_t$
is the smallest number of times the PDAE must be differentiated with respect to $t$
in order to determine $u_t$ as a continuous function of $t,~x,~u$ and certain
space derivatives of $u-$components.}
\end{defn}
Since in most practical applications of PDAEs the time index is $\nu_t=1$ or
$\nu_t=2$, we give here the formalism for these two indices for PDAEs
of type (\ref{pdae1}) only. First, an auxiliary result is stated. To simplify
the notation we use in the following the summation convention (summation about
twofold indices from $1$ to $n$).
\begin{lem} \label{hilf}
Let $u$ be sufficiently smooth, and suppose that  $C[u]$ is
linear in $u$. Then
\begin{enumerate}
\item [(I)] $\partial_tC[u]u_x=C^{(1)}[u_x]u_t+C[u]u_{tx}$
where
$C^{(1)}_{ik}[u_x]:=C_{ij,k}u_{j,x}$,
\item [(II)]
$\partial_tC^{(1)}[u_x]u_t=C^{(1)}[u_{tx}]u_t+C^{(1)}[u_x]u_{tt}.$
\end{enumerate}
\end{lem}
\begin{pf} Let $C_{ij,k}:=\frac{\partial C_{ij}}{\partial u_k}, ~i,j,k\in\{1,\dots,n\}.$
Componentwise differentiation with respect to time yields
\[\partial_tC_{ij}u_{j,x}=C_{ij,k}u_{k,t}u_{j,x}+C_{ij}u_{j,tx}=
(C_{ij,k}u_{j,x})u_{k,t}+C_{ij}u_{j,tx}.\]
Statement (I) follows from the definition of $C^{(1)}_{ik}$. Using this
definition again and the linearity of $C$ in $u$ we get
\begin{align*}
\partial_tC^{(1)}_{ik}[u_x]u_{k,t}=&
C_{ij,k}u_{j,tx}u_{k,t}+C_{ij,k}u_{j,x}u_{k,tt}
=&C^{(1)}_{ik}[u_{tx}]u_{k,t}+C_{ik}^{(1)}[u_x]u_{k,tt}
\end{align*}
which is (II) in component form.\qed
\end{pf}
For example, $C[u]$ given in (\ref{matrices}) produces
$C^{(1)}[u_x]=\tsm{u_{2,x}&u_{1,x}&0&0\\
0&u_{2,x}&0&0\\0&0&u_{4,x}&0\\0&0&0&0}.~$\\
With this Lemma and with definitions
\begin{align*}
L:=& B\partial^2_x+D,\qquad M[u,u_x]:=L+C^{(1)}[u_x]+C[u]\partial_x
\end{align*}
we find formally from eq. (\ref{pdae1}) under the assumptions of Lemma
\ref{hilf}
after one time differentiation (this case is relevant for $\nu_t=1$) the
derivative array
\begin{align} \label{index1}
\ma{A&0\\M[u,u_x]&A}\ma{u_t\\u_{tt}}+\ma{Lu+C[u]u_x\\0}
=&\ma{f\\f_t}.
\end{align}
Two time differentiations of (\ref{pdae1}) produce the derivative array
(for $\nu_t=2$)
\begin{align} \label{index2}
&\ma{A&0&0\\M[u,u_x]&A&0
\\2C^{(1)}[u_{tx}]&M[u,u_x]&A}
\ma{u_t\\u_{tt}\\u_{ttt}} +\ma{Lu+C[u]u_x\\0\\0}
=\ma{f\\f_t\\f_{tt}}.
\end{align}
When $A$ is singular, the coefficient matrices (denoted by $\mathcal{A}$)
of the  systems (\ref{index1}), (\ref{index2}) are
singular in the sense that ${\mathcal{A}}v=0$ where $v\neq 0$ is an appropriate
vector function (e.g. in the case of eq. (\ref{index1}), $v=(u_t^T, u_{tt}^T)^T,~u\in C^{2,2}_n$).
However, the coefficient matrices  might
uniquely determine $u_t$. The analogous problem for linear time
varying DAEs (differential algebraic
equations) is discussed in \cite[p. 29]{Bren89}, and for linear PDAEs of
first order it is investigated in \cite{Martin98}.\\
For the nonlinear PDAE (\ref{pdae1}), the derivative array (\ref{index1})
or (\ref{index2}) is used to write
\begin{align}\label{P1}
Pu_t=&F
\end{align}
where $P$ is according to eq. (\ref{index1}) or (\ref{index2}) an  operator
valued matrix,
and the vector $F$ does not depend on $u_t$.
First, we must find $P$ and an appropriate domain of definition $D(P)$. Second,
the invertibility of $P$ must be studied.\\
We mention, that the matrix $P$ in equation (\ref{P1}) is the analogue
to the nonsingular diagonal matrix $D(x_j)$ defined in \cite{Martin98}, Definition 3.5.
However, here we need not require a diagonal form of $P$.
\subsection{Spatial index} \label{spaceindex}
If $B$ is singular, we suppose that the PDAE is written in quasilinear form
(with respect to second space derivatives), i.e. we transform $B$
according to
$\overline B=S_0BS_1^{-1}=\tsm{I_{m}&0\\0&0}$
where $S_0,~S_1$ are constant regular $(n,n)-$matrices.
\begin{defn} \label{defspace}
{\it The spatial index $\nu_x$ of a system with a regular matrix $B$
is defined to be zero. When $B$ is singular,
$\nu_x$ is the smallest number of times the PDAE
\begin{align*}
\overline A\overline u_t+\overline B\overline u_{xx} +\overline C[S_1^{-1}\overline{u}]\overline
u_x+\overline D\overline u=\overline f(t,x)
\end{align*}
($\overline u=S_1u,~\overline A=S_0AS_1^{-1}$, and so on)
must be differentiated with respect to $x$ in order to obtain
\begin{align*}
\overline U:=(\underbrace{\overline u_{1,xx},
\overline u_{2,xx},...,\overline u_{m,xx}}_{m~~ elements},
\underbrace{\overline u_{{n_1+1},x},...,\overline u_{n,x}}_{n-m~~ elements})^T
\end{align*}
as a continuous function of $~~t,~x,~\overline u,~\overline u_t$ and
 $\bar{u}_{1,x},\dots,\bar{u}_{m,x}$.}
\end{defn}
It is straightforward to derive for
$\overline {U}$ a derivative array (by differentiations of the PDAE with
respect to $x$) analogous to eq. (\ref{index1}) or (\ref{index2}). As a result,
after a minimal number of $x-$differentiations one can find
a representation of $\overline{U}$ of the form
$\overline Q~ \overline U=\overline G$
 where $\overline Q$ is an operator valued matrix. The vector $\overline G$
is independent of the components of $\overline U$. This definition of the
spatial index is according to the one given in \cite{Lucht99}. It does not
transform the PDAE (\ref{pdae1}) to a system of first order as is often
done, e.g. in \cite{Martin98}. Such a transformation, in general,  changes
the indices. For example, one can show that the time-index here is $0$ if
and only if the corresponding index in \cite{Martin98} is $1$, a
time-index $1$ or $2$ here implies an index $2$ there.
\subsection{Time index of the plasma PDAE} \label{plindex}
Here, we study the time index $\nu_t$ of the PDAE
(\ref{pdae1}) under the assumption of zero BVs, $u(t,0)=u(t,1)=0,~$ and $f\in
C_{4}^{p,q}$ with suitable $p,~q.~$ To be specific we further assume that this
IBVP has a solution $u\in C_{4,0}^{1,2}$ which  possibly exists only local
in time. To determine the time index of system (\ref{pdae1}) we generate
for $\nu_t$  -- if possible -- system (\ref{P1}) with the condition that
$P$  defined on  $D(P)$ is invertible.
In order to determine $\nu_t$ we differentiate the third and fourth
equation of the PDAE with respect to $t$ with  the result
\[-u_{3,tx}+u_{3,t}u_{4,x}+u_3u_{4,tx}=0\;\;
\hbox{and}\;\;u_{4,txx}+u_{1,t}-u_{3,t}=0.\]
 A new system can be formed by
means of these two equations and the first two equations of the original
system. This is a closed system of four equations for the components of
$u_t$ in terms of $u$ and derivatives which are not time derivatives of
$u-$components. Obviously, the new system  can be written
\begin{align}\label{P3}
Pu_t=F(u):=&\tsm{-u_2u_{1,x}-u_1u_{2,x}+b_0u_{1,xx}\\
-u_2u_{2,x}-d_1u_{4,x}\\0\\0},
~P:=&\tsm{1&0&0&0\\
0&1&0&0\\0&0&[-\partial_x+u_{4,x}]&u_3\partial_x\\1&0&-1&\partial^2_x}.
\end{align}
Now what is essential  for the determination of the index of the PDAE is
that $u$ in $P$ can be seen as a fixed element with the result
that $P:~D(P)\rightarrow R(P)$, $D(P)\subset C_{4,0}^{1,2}$, is a  linear operator ($R(P)$
denotes the range of $P$).
We can try to solve the first equation in (\ref{P3}) for $u_t$ by
standard linear theory. In particular, if $P^{-1}$ exists we get $\nu_t$
as the least number of time differentiations which are necessary to get
the equation for $u_t$ (in the example considered here, one differentiation
with respect to $t$ is needed).
\begin{rem} The right hand side of the first equation in (\ref{P3})
 does not play
any role when $\nu_t$ is to be determined. This implies that $\nu_t$
is independent of $f$ in equation (\ref{pdae1}). This  is a
reasonable result  because
the right hand side function should not influence the indices of a PDAE.
\end{rem}
We ask whether the linear operator $P$ with $D(P)\subset C_{4,0}^{1,3}$
 has
for fixed $u=u(t,x)$ an inverse. The answer comes from  kernel $N(P)$ whose
elements $z$ fulfill
\begin{align*}
Pz=\tsm{z_1\\z_2\\-z_{3,x}+u_{4,x}z_3+u_3z_{4,x}\\
z_1-z_3+z_{4,xx}}=0.
\end{align*}
This  reduces to $z_1=z_2=0$ and a linear homogeneous
coupled system of two ordinary differential equations for
$z_3$ and $z_4,$
$-z_{3,x}+u_{4,x}z_3+u_3z_{4,x}=0$ and $-z_3+z_{4,xx}=0,$
with homogeneous BVs, $z_3(0)=0$ and $z_4(0)=z_4(1)=0.~$ $~z_4$ is also a solution
of the equation $-z_{4,xxx}+u_{4,x}z_{4,xx}+u_3z_{4,x}=0$ which can be reduced to
$-y_{xx}+u_{4,x}y_x+u_3y=0$ where $y=y(t,x):=z_{4,x}(t,x).$ By $\varphi_1$ and
$\varphi_2$ we denote two fundamental solutions of the equation for $y$, i.e.
$y=K_1\varphi_1+K_2\varphi_2$ is the general solution ($K_i,~i=1,2,$ are
independent of $x$). Then the following Lemma holds:
\begin{lem}\label{lemma1}
Let $u$ be such that
$\left| \sm{\varphi_{1,x}(0)&\varphi_{2,x}(0)\\
\int\limits_0^1\varphi_1(\xi)d\xi&\int\limits_0^1\varphi_2(\xi)d\xi}\right |~\neq~0.$\\
Then $z_3~=z_4~=~0.$ \qquad
\end{lem}
Since the proof is simple, it is omitted here (for details see \cite{Lucht001}).\\
Under the assumptions of this Lemma we see that $N(P)=\{0\}$. Therefore, system
(\ref{P3}) can be solved for $u_t$, and the time index is
$\nu_t=1$.\\
We mention that based on Definition \ref{defspace} the spatial index of the plasma PDAE can be determined
similarly. The result is $\nu_x=0.~~$
\section{Numerical solution by a  finite difference method} \label{splitting}
In this section we consider the numerical solution of IBVPs (\ref{pdae1}) --
(\ref{com}) by means of a fractional step difference method which
is combined with a matrix factorization. The fractional step method and  numerous
variants of it  are well known, see e.g. \cite{Jan69}. By this method, the order of
the system of equations  which must be solved can be reduced considerably.
The effort can be reduced further by a proper partition of the
original system matrix (denoted by $L$ below) into two splitting matrices
($L=L_1+L_2$).\\
To describe the general procedure, we
first rewrite  the nonlinear part $C[u]u_x$ of the PDAE (\ref{pdae1})
as $C[u,\partial_x]u$ which is more convenient sometimes.
$C[u,\partial_x]$ is an operator valued matrix. For example, with
$C[u]$ given in (\ref{matrices}) it is
\begin{align*}
C[u]u_x=\ma{u_2&u_1&0&0\\0&u_2&0&d_1\\0&0&-1&u_3\\0&0&0&0}\tma{u_{1,x}\\
u_{2,x}\\u_{3,x}\\u_{4,x}}=&\ma{\partial_xu_2&0&0&0\\
0&u_2\partial_x&0&d_1\partial_x\\0&0&-\partial_x&u_3\partial_x\\0&0&0&0}
\ma{u_1\\u_2\\u_3\\u_4}=:
C[u,\partial_x]u.
\end{align*}
We suppose that $A$ is singular and is given as
$A=\tsm{I_{n_1}&0\\0&0}$
~where $I_{n_1}$ is the unit matrix of order $n_1<n$ $(n_1\geq 1).$ Let
$n_2=n-n_1.$ Corresponding to this partition of $A$ we introduce the
notation $u=(u_1,~u_2)^T~$ and
$M=\tsm{M_{11}&M_{12}\\M_{21}&M_{22}}$ where the $(n,n)-$matrix $M$ may be one of
the matrices $B, ~C,~D$. The component representation of $u_k$ is
$u_k=(u_{k1}^T,...,u_{kn_k}^T)^T,~k=1,2.$
According to this partition,  PDAE  (\ref{pdae1}) is written
\begin{align}\label{pdae4}
\ma{I_{n_1}&0\\0&0}u_t+\ma{B_{11}&B_{12}\\B_{21}&B_{22}}u_{xx}
+\ma{C_{11}[u,\partial_x]&C_{12}[u,\partial_x]\\C_{21}[u,\partial_x]&
C_{22}[u,\partial_x]}u+\ma{D_{11}&D_{12}\\D_{21}&D_{22}}u=f.
\end{align}
Furthermore, let the operator valued matrices
$L_k[u],~k=1,2,~$ be defined by  ($C_{kj}=C_{kj}[u,\partial_x]$)
\begin{align}\label{L1}L_1[u]:=&\ma{0&0\\(B_{21}\partial^2_x+C_{21}+
D_{21})&(B_{22}\partial^2_x+C_{22}+D_{22})},\\ \label{L2}
L_2[u]:=&\ma{(B_{11}\partial_x^2+C_{11}+D_{11}) &
(B_{21}\partial^2_x+C_{21}+D_{21})\\0&0}
\end{align}
such that $L[u]~:=~B\partial^2_x+C[u,\partial_x]+D~=~L_1[u]+L_2[u].$
\begin{rem} Sometimes, other operators $L_1$, $L_2$ with $L=L_1+L_2$
may be more convenient  because the solution of the equations
may be easier. In every case, the partition should be such that the
identity $AL_2=L_2$ does hold. This is needed in a factorization as
explained below in eq. (\ref{factor}).
\end{rem}
In order to obtain an approximate numerical solution of IBVP (\ref{pdae1})
-- (\ref{com}) we consider it
with zero BVs (\ref{BC}) by means of a difference method.
The first time derivative is approximated with an equidistant time step size
$\tau$ by
$$u_t(t_{m+1},x_k)\approx\frac{1}{\tau}(u^{m+1}(x_k)-u^m(x_k)),\quad t_{m+1}=(m+1)\tau,
\quad m=0,1,...$$
The corresponding
IBVP is space discretized  on an equidistant grid
\begin{align*}
\Omega_h:=\{x_k=kh,~k=0,1,...,M;~~h=1/M,~M>1\}
\end{align*} by using difference formulas ($k=1,...,M-1$)
\begin{align}\nonumber
u_{xx}(t,x_k)\approx &\frac{\delta^2}{h^2}u_k(t):=
\frac{1}{h^2}\big(u_{k-1}(t)-2u_k(t)+u_{k+1}(t)\big),\\ \label{diff1}
u_x(t,x_k)\approx &\frac{\delta_0}{2h}u_k(t)\quad \hbox{or}\quad
\frac{\delta_+}{h}u_k(t)\quad\hbox{or}\quad\frac{\delta_-}{h}u_k(t)
\end{align}
where $\delta_0,~\delta_+$ and $\delta_-$ are the usual central, forward
and backward difference operators, respectively.
Now suppose  $t\in(0,t_e)$, $t_e>0$. We approximate
eq. (\ref{pdae1}) by the difference equation ($m=0,1,...,~~~k=1,...,M-1$)
\begin{align}\label{pdae5}
A\frac{u_k^{m+1}-u_k^m}{\tau}+L_h[u_k^m]
u_k^{m+1}=&f_k^{m+1}.
\end{align}
$L_h[u_k^m]$ denotes a discretization of $L[u(t_m,x_k)].$
We rewrite this equation  as
\begin{align}\label{pdae6}
\left(A+\tau
L_h[u^m_k]\right)\frac{u_k^{m+1}-u_k^m}{\tau}=f_k^{m+1}-L_h[u^m_k]u^m_k
\end{align}
and factorize approximately for $\tau \rightarrow 0$
\begin{align}\label{factor}
&A+\tau L_h[u_k^m]\approx (A+\tau
L_{h1}[u_k^m])(I+\tau
L_{h2}[u_k^m])\\ \nonumber
=&A+\tau(L_{h1}[u_k^m]+AL_{h2}[u^m_k])+O(\tau^2)
=A+\tau L_h[u^m_k]+O(\tau^2),
\end{align}
where $L_{h1}$ and $L_{h2}$ are corresponding discretizations of
(\ref{L1}) and (\ref{L2}), respectively (we used $AL_{h2}=L_{h2}$ and $L_{h1}+L_{h2}=L_h$).
$u_k^m$ is considered to be an approximation for $u(t_m,x_k)$.
Therefore, we study for $m=0,1,...$, $k=1,...,M-1$ fractional splitting
\begin{align}\label{spl1}
\left(A+\tau L_{h1}[u^m_k]\right)u_k^{m+1/2}=&f_k^{m+1}-L_h[u^m_k]u_k^m,\\
\label{spl2} \left(I+\tau
L_{h2}[u^m_k]\right)\frac{u_k^{m+1}-u_k^m}{\tau}=&u_k^{m+1/2}.
\end{align}
$u_k^0$ is (for every $k$) given as IV.
\begin{rem}\label{redorder}
This discussion shows that the fractional step method requires
the solution of two linear systems of coupled equations, but each of them is, in
numerous applications, of a considerably reduced order.
\end{rem}
The approximate factorization (\ref{factor}) implies the following Lemma.
\begin{lem}\label{lemmaspl}
Suppose
\begin{enumerate}
\item $v=v(t,x)\in{\Rset}^n,~t\in (0,t_e),$ is for some $t_e>0$
the exact solution (sufficiently smooth) of the IBVP (\ref{pdae4}),
(\ref{IC}) -- (\ref{com}),
\item the system of algebraic equations (\ref{spl1}), (\ref{spl2}) has a unique
solution.
\end{enumerate}
Then the method (\ref{spl1}), (\ref{spl2}) and the system (\ref{pdae6})
are equivalent for $\tau \rightarrow 0$, i.e. eqs. (\ref{spl1}), (\ref{spl2})
approximate the original system to the same $\tau-$order as eq. (\ref{pdae6}).
\end{lem}
\subsection{Convergence of the fractional step method}\label{full}
We consider the convergence of the numerical solution $u_k^m$ calculated by
the scheme (\ref{pdae5}) towards the exact solution
$v(t_m,x_k)$ of the corresponding IBVP for $\tau \rightarrow 0$ and
$h\rightarrow 0$ under the condition $~(m+1)\tau=t~~ (t~~\mbox{fixed}).~$
The basic assumption is that
$C[v]=C^0+C^1[v]$ is linear in $v$ where the matrix $C^0=const.$
 is chosen in such a manner
that the matrices $G_{0k}$, \mbox{$k=1,...,M-1,$} defined below in
(\ref{G0k}) are regular. For example, a proper choice may be
\mbox{$C^0=C[\bar v],$} $C^1[v]=C[v]-C[\bar v]$
where vector $\bar v$ does not depend on $t$ and $x,$ e.g. $\bar v$
is a suitable  mean
value of $v$.\\ [1ex]
First, we need the full truncation error $\alpha_k^{m+1}$ defined
for  $~m=0,1,...,$\\ $k=1,...,M-1$ by
\begin{align}\label{te1}
\alpha_k^{m+1}:=&A\frac{v_k^{m+1}-v_k^m}{\tau}+L_h[v_k^m]
v_k^{m+1}-f_k^{m+1}
\end{align}
where $L_h[v^m_k]:=B\frac{\delta^2}{h^2}+C[v_k^m]\frac{\delta}{qh}+D,$
and $q=1$ (for one--sided differences) or $q=2$ (for central
differences). $\delta$ stands for $\delta_0$ or $\delta_+$
or $\delta_-.$
Under the assumption that $v$ is
sufficiently smooth we Taylor expand $v^m_k,~v^m_{k\pm 1}$
in $t_{m+1}$ and $x_k$
to get in lowest order with respect to $\tau$ and $h$
\begin{align*}
\alpha_k^{m+1}=&A(v_{k,t}^{m+1}-\frac{1}{2}v_{k,tt}^{m+1}\tau+O(\tau^2))+B(v_{k,xx}^{m+1}+
\frac{1}{12}v_{k,xxxx}^{m+1}h^2+O(h^4))\\   &+C[v_k^{m+1}
-v_{k,t}^{m+1}\tau+O(\tau^2)](v_{k,x}^{m+1}+\bar O^{m+1}(h^q))
+Dv_k^{m+1}-f_k^{m+1}.
\end{align*}
Since $C[v_k^{m+1}-v_{k,t}^{m+1}\tau]=C^0+C^{1}[v_k^{m+1}]-\tau C^{1}[v_{k,t}^{m+1}]~$
and\\
$~Av_{k,t}^{m+1}+Bv_{k,xx}^{m+1}+\big(C^0+C^{1}[v_k^{m+1}]\big)v_{k,x}^{m+1}+Dv_k^{m+1}=f_k^{m+1},~$
the truncation error can be written
\begin{align}\label{te2}
\alpha_k^{m+1}=&A\bar O^{m+1}(\tau)+B\bar O^{m+1}(h^2)+
C[v_k^{m+1}]\bar O^{m+1}(h^q)\\ \nonumber
&+C^{1}[v_{k,t}^{m+1}]\bar O^{m+1}(\tau h^q)+C^{1}[v_{k,t}^{m+1}]\bar O^{m+1}(\tau).
\end{align}
Here, e.g., $\bar O^{m+1}(\tau)=\tau w,~w\in \Rset^n,~w$
independent of $\tau$.
Second, scheme (\ref{pdae5}) is written
by means of the Kronecker product in matrix form,
\begin{align} \label{full1}
\Big(I_{M-1}\otimes\frac{A}{\tau}+&Q_h[U^m]\Big)U^{m+1}=\Big(I_{M-1}\otimes\frac{A}{\tau}
\Big) U^m+F^{m+1},\\\nonumber
Q_h[U^m]:=&\frac{1}{h^2}P\otimes B+\frac{1}{qh}\tilde P_q\otimes
C[U^m]+I_{M-1}\otimes D
\end{align}
where $U^m:=(u_1^{mT},...,u_{M-1}^{mT})^T,\quad F^m:=(f_1^{mT},...,f_{M-1}^{mT})^T,~$  and
\begin{align} \label{P}
P:=&\tsm{-2 & 1 & &\\ 1& -2 & 1 &\\ & & \ddots &\\ & & 1 & -2 &}\in\Rset^{(M-1)\times (M-1)}.
\end{align}
$q$ and $~\tilde P_q~$ depend on the discretization of the first space derivative.
For example, using central difference formula in (\ref{diff1}), it is $q=2$ and
\[ \tilde{P}:=  \tsm{0  &1 & &\\ -1&  0 & 1&\\ & & \ddots &
\\ & & -1 &  0 &}\in\Rset^{(M-1)\times (M-1)}~\] (another differencing is chosen in (\ref{onesided})).
The $(n(M-1),n(M-1))-$matrix $\tilde P_q\otimes C[U^m]$ is a block matrix whose
block at position $(j,k)$ is $\tilde P_{q,jk}C[u_j^m],$ $u_j^m\in \Rset^n.$
Eqs. (\ref{te1}) and (\ref{te2}) imply that $~V^{m+1}$ satisfies equation
\begin{align} \label{full2}
\Big(I_{M-1}\otimes\frac{A}{\tau}+&Q_h[V^m]\Big)V^{m+1}=\Big(I_{M-1}\otimes\frac{A}{\tau}
\Big) V^m+F^{m+1} +{\mathcal{A} }^{m+1},\\
\label{calA}
{\mathcal{A}}^{m+1}:=&E_1O^{m+1}(\tau)
+E_2O^{m+1}(h^2)+E_3^{m+1}O^{m+1}(h^q)\\ \nonumber
&+E_4^{m+1}\left[O^{m+1}(\tau h^q) +O^{m+1}(\tau)\right]
\end{align}
where
$E_1:=I_{M-1}\otimes A,~~ E_2:=I_{M-1}\otimes B,~~
E_3^{j+1}:=I_{M-1}\otimes C[V^{j+1}],$\\ $ E_4^{j+1}:=I_{M-1}\otimes
C^1[V_t^{j+1}]$
and, e.g., $O^{m+1}(\tau)$ means $O^{m+1}(\tau)\in{\Rset}^{n(M-1)}.$
The foregoing relations can be used to estimate a norm
of the global error $\eta^m:=V^m-U^m$. In the following we choose the
discrete $L_2-$norm defined for a vector $v=(v_1,...,v_{n(M-1)})^T$ by
 $||v||:=\left[h\sum[k=1,n(M-1)]{v_k^2}\right]^{1/2}.~$\\
Subtracting  eq. (\ref{full1}) from eq. (\ref{full2}), we obtain
\begin{align*}
\left(I_{M-1}\otimes\frac{A}{\tau}\right)\eta^{m+1}+Q_h[V^m]V^{m+1}-Q_h[U^m]U^{m+1}=&\left(
I_{M-1}\otimes\frac{A}{\tau}\right)\eta^m+{\mathcal{A}}^{m+1}.
\end{align*}
By the identity
$$Q_h[V^m]V^{m+1}-Q_h[U^m]U^{m+1}=\big(Q_h[V^m]-Q_h[U^m]\big)
V^{m+1}+Q_h[U^m]\big(V^{m+1}-U^{m+1}\big),$$
the foregoing equation can be written
\begin{align}\label{et1}
G^m\eta^{m+1}=&\left(I_{M-1}\otimes\frac{A}{\tau}\right)\eta^m-\left(Q_h[V^m]
-Q_h[U^m]\right)V^{m+1}+{\mathcal{A}}^{m+1},\\  \label{G} G^m:=&I_{M-1}\otimes\frac{A}{\tau}
+\frac{1}{h^2}P\otimes B+\frac{1}{qh}\tilde P_q\otimes C[U^m]+I_{M-1}\otimes D.
\end{align}
Now we use the fact that \[Q_h[V^m]-Q_h[U^m]=\frac{1}{qh}\tilde
P_q\otimes\left(C[V^m]-C[U^m]\right)=\frac{1}{qh}\tilde
P_q\otimes C^{1}[\eta^m]\] (because the linearity of $C[u]$  implies
$C[V^m]-C[U^m]=C^{1}[\eta^m]$). Furthermore, one can
construct a $(n(M-1),n(M-1))-$matrix $\tilde C[V^{m+1}]$ such that
\[\left(\tilde P_q\otimes C^{1}[\eta^m]\right) V^{m+1}=\tilde C[V^{m+1}]\eta^m.\]
Therefore, eq. (\ref{et1}) takes the form ($G^{-m}=(G^m)^{-1}$)
\begin{align*}
\eta^{m+1}=& G^{-m}\left(I_{M-1}\otimes\frac{A}{\tau} -
\frac{1}{qh}\tilde C[V^{m+1}]\right)\eta^m+G^{-m}{\mathcal{A}}^{m+1}
\end{align*}
or for short
\[~\eta^{m+1}=H^m\eta^m+R^{m+1},\quad
H^m:=G^{-m}\left(I_{M-1}\otimes\frac{A}{\tau}-\frac{1}{qh}\tilde
C[V^{m+1}]\right),\]\\ $\quad R^{m+1}:=G^{-m}{\mathcal{A}}^{m+1}.~$  It
follows for $m=0,1,...$
$$\eta^{m+1}=H^mH^{m-1}...H^0\eta^0+
H^mH^{m-1}...H^1R^1
+...
+H^mR^m+R^{m+1},$$
\begin{align}\label{et3}
||\eta^{m+1}||\leq &||H^mH^{m-1}...H^0||~||\eta^0||+
||H^mH^{m-1}...H^1||~||R^1||\\ \nonumber &+||H^mH^{m-1}...H^2||~||R^2||+...
+||H^m||~||R^m||+||R^{m+1}||.
\end{align}
Now we require for  $0<m\tau=t\in(0,t_e]$ that
\begin{align}\label{Lax}
\sup_{m\in{\mathbb{N}}}\left\{||H^mH^{m-1}...H^j||,~~j=1,...,m
\right\}~<~\infty,
\end{align} possibly under a restriction of one of the forms
\begin{align}\label{restrict}
\kappa_0\leq\frac{\tau}{h}\leq\kappa_1,&\qquad
\kappa_2\leq \frac{\tau}{h^2}\leq\kappa_3.
\end{align}
Assumption (\ref{Lax}) is discussed shortly in
 Remark \ref{lb}. With this condition,   the
right hand side in (\ref{et3}) can be estimated further to give
\begin{align}\nonumber
||\eta^{m+1}||&\leq K_0||\eta^0||+\bar K_1\sum[j=0,m]{||R^{j+1}||}\leq
K_0||\eta^0||+(m+1)\bar K_1\max_{j\in[0,m]}||R^{j+1}||\\ \label{et4}
&\leq K_0||\eta^0||+K_1\frac{t}{\tau}\max_{j\in[0,m]}||R^{j+1}||.
\end{align}
$K_0,~K_1,~\bar K_1$ are constants. Using
$||R^{j+1}||=||G^{-j}{\mathcal{A}}^{j+1}||,$ relation  (\ref{calA}) yields
\begin{align} \label{R0}
||R^{j+1}||\leq& ||G^{-j}E_1||O(\tau)+||G^{-j}E_2||O(h^2)+
||G^{-j}E_3^{j+1}||O(h^q)\\\nonumber
&+||G^{-j}E_4^{j+1}||\left[O(\tau h^q)+O(\tau)\right]
\end{align}
or also
\begin{align} \label{R1}
||R^{j+1}||\leq&||G^{-j}||\Big(||E_1||O(\tau)+||E_2||O(h^2)+
||E_3^{j+1}||O(h^q)\\\nonumber&+||E_4^{j+1}||\left[O(\tau h^q)+O(\tau)\right]\Big).
\end{align}
Note that ineq. (\ref{R0}) is sharper than ineq. (\ref{R1}).
Since the matrices $E_l$ are block diagonal, their norms can be estimated easily. For
example, $E_3^{j+1}$ is of the form\\ $E_3^{j+1}=diag\{C[v_1^{j+1}],...,C[v_{M-1}^{j+1}]\}$, and
$$||E_3^{j+1}||=\left[
\max_k\left[ \lambda_{max}(C^T[v_k^{j+1}]C[v_k^{j+1}])\right]\right]^{1/2}$$
where $\lambda_{max}(E)$ is the largest eigenvalue of $E$. Especially,
if $v$ is sufficiently smooth for some time $t\leq t_e$, then all norms $||E_l||$
are bounded. Therefore, (\ref{R1}) can be written also
\begin{align*}
||R^{j+1}||\leq K_3||G^{-j}||\left(O(\tau)+O(h^2)+O(h^q)+O(\tau h^q)\right)
\end{align*} with some constant
$K_3$.\\
An  estimate of $||G^{-j}||$ can be  based on the eigenvectors
$\phi_k$ of the symmetric matrix $\frac{1}{h^2}P$. By definition,
$\frac{1}{h^2}P\phi_k=\lambda_k\phi_k,~~k=1,...,M-1$.  The eigenvalues are
given by $\lambda_k=-\frac{4}{h^2}\sin^2\left(\frac{k\pi}{2M}\right)$.
They fulfill for $h\to 0$ the asymptotic relation
$\lambda_k=-(k\pi)^2+O(h^2)$ (because $h=1/M$). The set
$\{\phi_k,~~k=1,...,M-1\}$ is assumed to be orthonormal. Let
$\Phi:=\sqrt{h}(\phi_1,...,\phi_{M-1})$ be the matrix of the eigenvectors
and $\Lambda:=diag(\lambda_1,...,\lambda_{M-1}).$ Since
$\Phi^{-1}=\Phi^T=\Phi,~$ it follows that
(for short $\Psi_n:=\Phi\otimes I_n$)
\begin{align}\label{PB}
\frac{1}{h^2}P\otimes B=&\Psi_n(\Lambda\otimes B)\Psi_n,\\
I_{M-1}\otimes\left(\frac{1}{\tau}A+D\right)\label{PB1}
=&\Psi_n
\left(I_{M-1}\otimes\left(\frac{1}{\tau}A+D\right)\right)\Psi_n.
\end{align}
For simplicity, we  assume that the first order derivative $\partial_x$ is
discretized by a one sided difference approximation ($q=1$), say
\begin{align}\label{onesided}
\frac{1}{h}\tilde P_q\equiv & \frac{1}{h}\tilde{P}=\frac{1}{h} \tsm{-1  &1
& &\\ &  -1 & 1&\\ & & \ddots &\\& & & -1&1\\& & & & -1&}\equiv
\frac{1}{h}\left(-I_{M-1}+H_{M-1}\right).
\end{align}
Inserting this into eq. (\ref{G}) and using eqs.
(\ref{PB}) and (\ref{PB1}), $G^j$ can be written
\begin{align}\label{Gj}
G^j=&G_0+G_1^j,\\G_0\label{G0j}
 =&\Psi_n\left[I_{M-1}\otimes\left(\frac{1}{\tau}A
-\frac{1}{h}C^{0}+D\right)+\Lambda\otimes B\right]\Psi_n,\\
\label{G1j} G_1^j=&\frac{1}{h}\left[H_{M-1}\otimes C^{0}
+\tilde P\otimes C^1[U^j]\right].
\end{align}
Eq. (\ref{G0j}) implies that the $(n(M-1),n(M-1))-$matrix $G_0$ can be
reduced to the set of $(n,n)-$matrices $G_{0k},~k=1,...,M-1,$
\begin{align}
G_0=\Psi_n\tsm{G_{01}& & &\\& & \ddots &\\& & &
G_{0M-1}} \Psi_n,
 \label{G0k}
~~G_{0k}:=\frac{1}{\tau}A-\frac{1}{h}C^{0}+D+\lambda_kB.
\end{align}
If the low order
matrices $G_{0k},~k=1,...,M-1,$ are regular (this is due to the choice of $C^0$),
  then $G_0$ is regular. Therefore,
$G^j=G_0\left(I_{n(M-1)}+G_0^{-1}G_1^j\right),$ and it follows
\begin{align}\label{nGmj}
||G^{-j}||~\leq~
&||\left(I_{n(M-1)}+G_0^{-1}G_1^j\right)^{-1}||~||G_0^{-1}||
\end{align}
provided $I_{n(M-1)}+G^{-1}_0G^j_1$ is also invertible. The second
factor on the right hand side  can be written
\begin{align} \label{nG0mj}
||G_0^{-1}||=\max_k||G^{-1}_{0k}||_n
&~=~\max_k\frac{1}{\left[\lambda_{min}
\left(G_{0k}^{T}G_{0k}\right)\right]^{1/2}}
\end{align}
where $||Q||_n$ is the spectral norm of a real $(n,n)-$matrix $Q,$  and
$~\lambda_{min}\left(G_{0k}^{T}G_{0k}\right)$ is the smallest
eigenvalue of the matrix $G_{0k}^{T}G_{0k}.$ The first factor on the right of
ineq. (\ref{nGmj}) can be estimated by means of
\begin{align*}
\left|\left|\left(I_{n(M-1)}+G_0^{-1}G_1^j\right)^{-1}\right|\right| \leq
\frac{1}{1-||G_0^{-1}G_1^j||}
\end{align*}
provided $||G_0^{-1}G_1^j||<1.$ It may be
appropriate to estimate further\\ $||G_0^{-1}G_1^j||\leq
||G_0^{-1}||~||G_1^j||.~$ $~||G_0^{-1}||$ is given already in
(\ref{nG0mj}), and $||G_1^j||$ is
\begin{align}\nonumber
||G_1^j||
=&\frac{1}{h}\left|\left|-I_{M-1}\otimes C^1[U^j]+H_{M-1}\otimes
C[U^j]\right|\right|\\\nonumber
\leq& \frac{1}{h}\left(\left|\left|I_{M-1}\otimes C^1[U^j]\right|\right|+
\left|\left|H_{M-1}\otimes C[U^j]\right|\right|\right)\\ \label{CC}
=&\frac{1}{h}\left(\max_{k\in[1,M-1]}\left|\left|C^1[u_k^j]\right|\right|_n+
\max_{k\in[1,M-2]}\left|\left|C[u_k^j]\right|\right|_n\right).
\end{align}
Here we used the identities
\begin{align*}
&\left|\left|  H_{M-1} \otimes C[U^j] \right|\right| =
\left[\lambda_{max}\left(\left(H_{M-1}\otimes C[U^j]\right)^T
\left(H_{M-1}\otimes C[U^j]\right)\right)\right]^{1/2}\\
&=\left[~\lambda_{max}\tsm{0 & & & & \\& &  C_1^TC_1 &
& \\& &  &\ddots & &\\& & &
&C^T_{M-2}C_{M-2}}~\right]^{1/2}=
\max_k\left[~\lambda_{max}\left(C_k^TC_k\right)\right]^{1/2}
=\max_k||C_k||_n
\end{align*}
(for short $C_k=C[u_k^j]$).  Relations (\ref{nG0mj}) and (\ref{CC}) show that
the norms $||G_0^{-1}||$ and $||G_1^j||$ can be estimated by norms of low order
matrices (of order $n$ where $n$ is usually small, e.g. $n=4$ for the
plasma PDAE given in section \ref{plasma}).\\
We summarize the foregoing result in the following Lemma.
\begin{lem}\label{CFL1}
Suppose that for $\tau\in(0,\tau_0]$ and $h\in(0,h_0]$
($\tau_0,~h_0>0$),\\
$(m+1)\tau=t\in(0,t_e],$
\begin{enumerate}
\item the $(n,n)-$matrices $G_{0k},~k=1,...,M-1,$ are regular,
\item $||u_k^j||_n\leq K_0,~k=1,...,M-1,~j=1,...,m+1,$ where $K_0$ is some
constant independent of $\tau$ and $h,$
\item there is a positive constant $\delta_0<1$ such that the condition
\begin{align}\label{CFL2}
||G_0^{-1}||~||G_1^j||\leq&\Big(\max_{k \in[1,M-1]}~||G_{0k}^{-1}||_n\Big)~
\frac{1}{h}\Big(\max_{k\in[1,M-1]}~||C^1[u_k^j]||_n\\ \nonumber
&+\max_{k\in[1,M-2]}\left|\left|C[u_k^j]\right|\right|_n\Big)
~\leq \delta_0
\end{align}
is satisfied.
\end{enumerate}
Then
$\left|\left|G^{-j}\right|\right|~\leq~\frac{1}{1-\delta_0}
\left|\left|G_0^{-1}\right|\right|.\qquad $
\end{lem}
This condition plays a crucial role when proving the  convergence of the
dif\-fe\-rence scheme considered here. We illustrate this by the following  example.
\begin{exmp}\label {CFL} Inequality (\ref{CFL2}) is for the case
$n=1$, $A=1$, $B=-1$, $C=C^0=const.<0,$ $D=0$ and Dirichlet BVs of type of a CFL--condition
known from the discretization theory of hyperbolic differential equations. In
this case
\begin{align*}
G&=\frac{1}{\tau}I_{M-1}-\frac{1}{h^2}P+\frac{C^0}{h}\tilde P=G_0+G_1,\\
G_0&=\Phi\left[\left(\frac{1}{\tau}-\frac{C^0}{h}\right)I_{M-1}
-\Lambda\right]\Phi,\quad G_1=\frac {C^0}{h}H_{M-1}.
\end{align*}
Therefore, $G_{0k}=\frac{1}{\tau}-\frac{C^0}{h}-\lambda_k$ which implies
$$||G_0^{-1}||=\frac{\tau}{\min\limits_k\left(1+|C^0|\frac{\tau}{h}+
\tau|\lambda_k|\right)}\leq \tau.$$
Furthermore, $H^T_{M-1}H_{M-1}$ is a diagonal matrix which has
eigenvalues $0$ and $1$, hence $||H_{M-1}||=1,$ and
$$||G_0^{-1}||~||G_1||\leq \tau\frac{|C^0|}{h}.$$
If this is required to be less than $1$ (according to (\ref{CFL2})),
it resembles the well known CFL condition $|C^0|\tau/h<1.$
\end{exmp}
In order to prove convergence, the
estimate (\ref{R0}) should be applied to the error
inequality (\ref{et4}). To estimate the norms $||G^{-j}E_l||,~l=1,...,4,$
we use again the decomposition (\ref{Gj}) under the assumption that
$G_0^{-1}$ exists. Then the relation
$G^j=G_0\left(I_{n(M-1)}+G_0^{-1}G_1^j\right)$ gives
\begin{align*}
\left|\left|G^{-j}E_l\right|\right|~=&
\left|\left|\left(I_{n(M-1)}+G_0^{-1}G_1^j\right)^{-1}
G_0^{-1}E_l\right|\right|\\
~\leq &~\left|\left|\left(I_{n(M-1)}+G_0^{-1}G_1^j
\right)^{-1}\right|\right|~\left|\left|G_0^{-1}E_l\right|\right|.
\end{align*}
Since $E_l$ is (for all $l$) a block diagonal matrix, the representation
(\ref{G0j}) implies that also the norm $\left|\left|G^{-j}E_l\right|\right|$
can be estimated by the calculation of norms of $(n,n)-$matrices,
\begin{align*}
||G^{-j}E_l||
\leq~&\left|\left|\left(I_{n(M-1)}+G_0^{-1}G_1^j
\right)^{-1}\right|\right|~\max_k~||G_{0k}^{-1}||_n\max_{k'}~||E_{lk'}||_n.
\end{align*}
Sometimes, this estimate may be useful.
The result of the foregoing estimates is the following Lemma.
\begin{lem}\label{et5}
 Suppose\begin{enumerate}
\item $C[u]=C^0+C^1[u]\in{\Rset}^{n\times n}~$ is linear in $u$,
\item the exact solution  $v=v(t,x)$ of the IBVP is sufficiently smooth for\\  $~t\in(0,t_e],~$
$t_e>0,~$ especially $||V||,~~||V_t||<\infty,$
\item $(m+1)\tau=t\in(0,t_e],~m\in\mathbb{N},~$ $t$ fixed,
\item $||U^j||<\infty,~j=1,...m,$
\item $||\eta^0||=0,$
\item $\sup_{m\in{\mathbb{N}}}\left\{||H^mH^{m-1}...H^j||,~~j=1,...,m
\right\}~<~\infty,$
\item inequality (\ref{CFL2}) with $\delta_0<1$ is valid.
\end{enumerate}
Then an upper bound  of the global error  (see ineq. (\ref{et4}))  can be
expressed in terms of norms containing $G_0^{-1}$ only:
\begin{align*}
||\eta^{m+1}||~\leq~&\frac{K_1}{1-\delta_0}~\frac{t}{\tau}\Big\{
\max_{k\in[1,M-1]}\Big[
 ||G^{-1}_{0k}A||_n ~O(\tau)
 \nonumber +||G^{-1}_{0k}B||_n~O(h^2)\Big]\\
 &+\max_{j\in[1,m+1]}\Big( ||G_0^{-1}E^{j}_3||O(h)
  +||G_0^{-1}E^{j}_4|| (O(\tau)+O(\tau h))\Big)
\Big\}.\qquad
\end{align*}
\end{lem}
From this Lemma, one may get convergence results possibly under a
restriction of the type given in (\ref{restrict}).
$||G_{0k}A||_n,~||G_{0k}^{-1}B||_n~$ and $||G_0^{-1}E_l||,~l=3,4,~$ depend on
$\tau,~h$ and on the indices of the PDAE. For example, in the case of the
plasma system ($\nu_t=1$), we have with the choice
$C^0=\tsm{0&0&0&0\\0&0&0&d_1\\0&0&-1&0\\0&0&0&0}$
$||G_{0k}^{-1}A||_n=O(\tau),~~||G_{0k}^{-1}B||_n=O(\tau),$ provided $\tau$
and $h^2$ are related by $\tau/h^2\leq \kappa$ where $\kappa$
is independent of $\tau$ and $h$. By the way we note that in the plasma example  $||G_{0k}^{-1}||_n$
is of order $O(\tau^{1/2})$ only. Therefore, one should not estimate
$||G_{0k}^{-1}A||_n$ by the upper bound $||G_{0k}^{-1}||_n||A||_n.~$ The dependence of
norms on the indices of the PDAE is discussed in \cite{Lucht99}.
\begin{rem}\label{lb}~
\begin{enumerate}
\item  The sixth requirement is a stability condition.
Conditions of this type are  well known in the theory of difference methods
of time dependent partial differential equations (see also \cite{Lucht99}
where the assumption is discussed for linear PDAEs without convection).
Unfortunately, a general easy method to verify the assumption when a
convection term is present  is not avai\-lable. It should be
studied separately for each problem of type (\ref{pdae1}) with given
matrices $A,~B,~C,~D.$
\item In most cases, the seventh assumption can be satisfied only when the
step sizes $\tau$ and $h$ are restricted in a certain way and/or when
the convection term is not dominant, see also example \ref{CFL}.
\end{enumerate}\end{rem}
\section{Numerical example} \label{numex}
For an example we choose the plasma PDAE
(with $n=4,~d=1$) described in section \ref{plasma}. The matrices are defined
in (\ref{matrices}).  Let the parameters in $B$ and $C[u]$ be defined by
$b_0=0.02$ and $d_1=1.$ First, according to (\ref{IC}) and (\ref{BC}) IVs and
(Dirichlet) BVs must be specified. Since $\nu_t=1$ (see section \ref{plindex}),
both terms $\Phi_a,~\Phi_c$ are different from zero. One can prescribe arbitrary
IVs for two components of $u=u(t,x)$ only, say $u_2$ and $u_4$ (these are the
non zero components of $\Phi_a$). One can show that we can assign especially
 the  BV $u_3(0,0)$ which is assumed
to be non zero. Given this BV and
\begin{align*}
g_4(x):=u_4(0,x)=K_4\cos(2\pi x)
\end{align*}
(where the constant $K_4$ is defined below), it follows
from the third equation of the system (\ref{pdae1}), (\ref{matrices}) that the
consistent IV of $u_3$ is
\begin{align*}
g_3(x):=u_3(0,x)=u_3(0,0)e^{g_4(x)}/e^{g_4(0)}.\end{align*}
This can be inserted into the fourth equation of the system with the result that
the consistent IV for $u_1$ is
\begin{align*}
g_1(x):=u_1(0,x)=g_3(x)-g_{4,xx}(x)=u_3(0,0)e^{g_4(x)}/e^{g_4(0)}+4\pi^2g_4(x).
\end{align*}
The IV for $u_2$ can be given arbitrarily, here we choose
\begin{align*}
g_2(x):=u_2(0,x)=K_2x(x-0.5),\qquad K_2=const.,
\end{align*}
and BVs are for $t\in \bar I_t$
\begin{align*}
u_1(t,0)&=u_1(t,1)=u_2(t,0)=0,~ u_3(t,0)=u_3(0,0),
~ u_4(t,0)=u_4(t,1)=K_4.
\end{align*}
According to relation (\ref{com}), the BVs for $u_1$ imply
 \mbox{$K_4=-u_3(0,0)/(4\pi^2).$} With these relations  the IBVP
(\ref{pdae1}) -- (\ref{com}) with matrices (\ref{matrices}) is defined
completely, and we may solve it by a finite difference method. The two
parameters $K_2$ and $u_3(0,0)$ can be chosen arbitrarily, for an example
let $u_3(0,0)=0.2$ and \mbox{$K_2=0.4$}. Furthermore,  let the two step sizes
$\tau$ and $h$ be related by \mbox{$\tau=K_0h$} (or less restrictive by
$\tau\leq K_0h$) where $K_0$ is a positive constant ($K_0=0.5$ below).
This restriction has its origin in the second equation which is of
hyperbolic type (inhomogeneous Burgers equation). In the example, the CFL
condition for the second equation is satisfied for the data chosen because
$$\frac{\tau}{h}\max_{k}\left\{|u^{m}_{2k}|_{m\tau=1}\right\}$$ is much
smaller than 1. Some values of this expression are given in Table 1
\begin{table}[ht]
\center
\caption{}
\begin{tabular}{|c*{5}{|c}|}\hline
N & 20 & 40 & 80 & 160 & 320\\\hline
CFL$_2$ & 0.0771 & 0.0799 & 0.0811 & 0.0817 & 0.0819\\
 e$_1$&0.0075&0.0051 & 0.0027 &   0.0013 &   0.0006\\
 e$_2$& 0.0141  &  0.0113 &   0.0076 &   0.0049 &   0.0031
\\\hline
\end{tabular}
\end{table}
in
the line CFL$_2$ for different values of $N=1/h$. Furthermore, in the
table $e_i$ is defined by $~e_i:=||U_{i,h}-U_{i,h/2}||,~i=1,2,~$ where
$U_{i,h}$ and $U_{i,h/2}$ are the numerical approximations of $u_i$ at
time $t=1$ and at the space grid points for the two different step sizes
$h$ and $h/2$. Note that $u_1$ and $u_2$ are the ion density and the ion
velocity, respectively (see section \ref{plasma}), and,  therefore, $u_1$
must be nonnegative.\\
In Figure 5.1, a typical result of the finite
difference solution of the IBVP with  the data given at time $t=1$ is
shown. The components $u_1$  and $u_2$  which are the most interesting
ones are presented.
\begin{figure}[ht]
\caption{}
    \epsfig{bbllx=13,bblly=13,bburx=599,bbury=599,clip=,figure=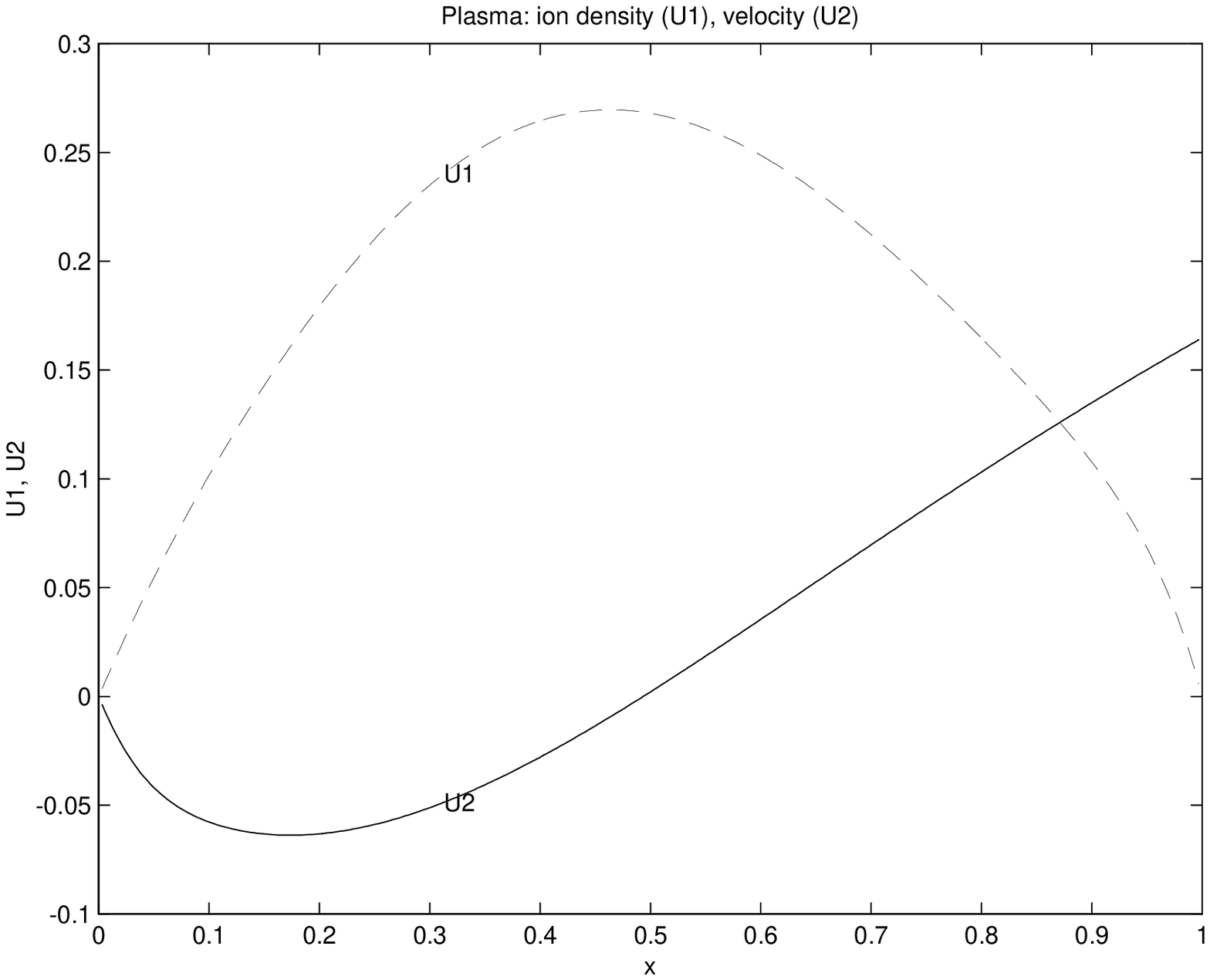,height=12cm
    }
\vspace{-3cm}
\end{figure}
\section{Conclusion}
After an application of PDAEs with a nonlinear convection term we first
consi\-dered the determination of the time and  spatial index of such systems. These were
based on definitions given earlier in the literature,  e.g. \cite{Martin98}.  We reduced this problem to the question whether
there is an inverse of some operator defined on a properly defined solution space.
 For the case of the plasma
PDAE (a system with time index 1) this was shown in detail.  In most practical applications,  PDAEs
have time index less than 3.\\ Then we studied the
numerical solution of corresponding IBVPs by means of a  finite difference splitting
method familiar from the treatment of
partial differential equations. Especially the convergence
was considered for the case that both the time and space step sizes tend to zero.
The main problem in these investigations is (in one space dimension) to handle
the  limit $h\to 0$ in relation to the time step size $\tau$. For fixed $h$
and $\tau \to 0$ this problem is easy. Some results concerning the problem when
both $h$ and $\tau$ go to zero were given.\\
Finally, some results of  a numerical solution of an IBVP  from plasma physics were
presented. This example is interesting because the PDAE which is nonlinear is a mixture of
partial differential equations of parabolic, elliptic and hyperbolic type.
\bibliographystyle{plain}

\end{document}